\documentclass[reqno, 12pt, draft]{amsart}
\usepackage{amsmath,amsthm,amsfonts,amssymb}
\usepackage[cp1251]{inputenc}
   \usepackage[english]{babel}
\date{}
\newtheorem{theorem}{Theorem}[section]
\newtheorem{lemma}[theorem]{Lemma}%[section]
\newtheorem{example}[theorem]{Example}%[section]
\numberwithin{equation}{section}

\title[Estimates for the distances between solutions]
{Estimates for the distances between solutions to Kolmogorov equations with diffusion matrices of low regularity}

\author[Vladimir I. Bogachev and Stanislav V. Shaposhnikov]{Vladimir I. Bogachev$^{1,2}$ and Stanislav V. Shaposhnikov$^{1,2}$}



\thanks{This research is supported by the
Russian Science Foundation  Grant 25-11-00007.}

\begin{document}

\maketitle

{\footnotesize
 \centerline{$^1$Department of Mechanics and Mathematics, Lomonosov Moscow State
University, Moscow, Russia}
}

{\footnotesize
  \centerline{$^2$National Research University Higher School of Economics, Moscow, Russia}
}

\begin{abstract} We obtain estimates for the weighted $L^1$-norm of the difference of two probability solutions to  Kolmogorov equations
in terms of the difference of the  diffusion matrices and the drifts. Unlike the previously known results, our estimate does not involve
Sobolev derivatives of solutions and coefficients. The diffusion matrices are supposed to be non-singular, bounded
and satisfy the Dini mean oscillation condition.
\end{abstract}

\section{Introduction}

In the papers \cite{BKSH14}, \cite{BKSH17}, and \cite{BRSH18}, some estimates were
 obtained for the distances  between probability solutions $\varrho_{\mu}$ and $\varrho_{\sigma}$
to two Kolmogorov equations (also called stationary Fokker--Planck--Kolmogorov equations)
$$
\partial_{x_i}\partial_{x_j}\bigl(a^{ij}_{\mu}\varrho_{\mu}\bigr)-\partial_{x_i}\bigl(b^i_{\mu}\varrho_{\mu}\bigr)=0, \quad
\partial_{x_i}\partial_{x_j}\bigl(a^{ij}_{\sigma}\varrho_{\sigma}\bigr)-\partial_{x_i}\bigl(b^i_{\sigma}\varrho_{\sigma}\bigr)=0.
$$
Here the usual convention about summation over repeated indices is employed.

For example, in \cite[Theorem 3.2]{BRSH18}, the following  estimate was obtained:
$$
\|\varrho_{\mu}-\varrho_{\sigma}\|_{L^1(\mathbb{R}^d)}\le
C\int_{\mathbb{R}^d}\bigl(1+|x|\bigr)^m|\Phi(x)|\varrho_{\sigma}(x)\,dx,
$$
where
$$
\Phi=\bigl(\Phi^i\bigr)_{1\le i\le d}, \quad
\Phi^i=\sum_{j=1}^d(a^{ij}_{\mu}-a^{ij}_{\sigma})\frac{\partial_{x_j}\varrho_{\sigma}}{\varrho_{\sigma}}+
\sum_{j=1}^d\bigl(\partial_{x_j}a^{ij}_{\mu}-\partial_{x_j}a^{ij}_{\sigma}\bigr)
+\bigl(b_{\sigma}^i-b_{\mu}^i\bigr).
$$
Let us  also mention the paper \cite{FM22}, in which an estimate for the Kantorovich distance between stationary distributions
was obtained under the assumption that the coefficients belong to the Sobolev class with respect to
one of the two measures, these measures are equivalent and
 the one-dimensional distributions of the cor\-res\-pon\-ding random process converge exponentially fast to the stationary
distribution in the Kantorovich metric. Related results for  parabolic equations have been recently obtained in \cite{HRW} and~\cite{RW}.

Suppose now that the elements $a^{ij}_{\mu}$ and $a^{ij}_{\sigma}$ of the
diffusion matrices $A_{\mu}$ and $A_{\sigma}$ have a modulus of continuity
satisfying  the Dini mean oscillation condition (recalled below).
Then even in the one-dimensional case solutions to the Kolmogorov equation can fail to have
Sobolev derivatives. For example, any solution to the equation
$(a\varrho)''=0$, where a continuous function $a$ is positive,
is given by the formula $\varrho(x)=(c_1+c_2x)/a(x)$. If the function $a$ is H\"older continuous, but nowhere differentiable,
then for $c_1=1$ and $c_2=0$ the solution is nowhere  differentiable, hence has no Sobolev derivative.
In this case the known results, in particular, the aforementioned estimate for the $L^1$-norm of the difference of solutions,
are not applicable.
In this paper, in the case of a continuous diffusion matrix satisfying the Dini mean oscillation condition
we obtain the estimate
\begin{multline*}
\|(1+|x|^k)(\varrho_{\mu}-\varrho_{\sigma})\|_{L^1(\mathbb{R}^d)}\le
C\|A_{\mu}-A_{\sigma}\|_{L^r(\mathbb{R}^d, \varrho_{\sigma}\,dx)}
\\
+
C\int_{\mathbb{R}^d}|b_{\mu}(x)-b_{\sigma}(x)|(1+|x|^{\beta+k})\varrho_{\sigma}(x)\,dx,
\end{multline*}
where $C$ is a constant depending on certain quantities appearing in our conditions on the coefficients. 
The derivation of this bound is based on several recent results from the papers
\cite{DKE1}, \cite{DKE2}, 
and \cite{BRSH23}, devoted to  Kolmogorov equations with
 continuous  diffusion matrices without Sobolev derivatives.
In addition, as an auxiliary object of independent interest  we
study the   properties of solutions to the  Poisson equation. Such equations with unbounded  coefficients
were investigated in \cite{Veret} and \cite{BRSH18}. In this paper, we extend the results  from \cite{BRSH18}
to the case where the  elements of the  diffusion matrix are continuous and possess a modulus of continuity satisfying 
the Dini mean oscillation condition.
Note that estimates for the distances  between solutions are useful for the study of dependence of solutions on parameters
and provide broad sufficient conditions for the existence and uniqueness of stationary  solutions
to nonlinear Fokker--Planck--Kolmogorov equations and nonlinear Mckean--Vlasov equations (see \cite{BKSH14} and \cite{BS24}).
For other recent applications of results on the Poisson equation and estimates for the distances  between solutions,
see \cite{BGNS} and~\cite{Kouhkouh}.
Surveys of results on Kolmogorov equations are given in  \cite{BKRS} and~\cite{BS25}.

\section{Notation and terminology}

The open ball of radius $r$ centered at $x$ is denoted by $B(x, r)$.
Following \cite{DKE1} and \cite{DKE2}, we shall say that a measurable function $f$ on $\mathbb{R}^d$
satisfies the Dini mean oscillation condition if for some $t_0>0$
\begin{equation}\label{din1}
\int_0^{t_0}\frac{\omega(t)}{t}\,dt<\infty,
\end{equation}
where
$$
\omega(r)=\sup_{x\in\mathbb{R}^d}\frac{1}{|B(x, r)|}\int_{B(x, r)}|f(y)-f_{B}(x, r)|\,dy, \
f_{B}(x, r)=\frac{1}{|B(x, r)|}\int_{B(x, r)}f(y)\,dy,
$$
and $|B(x, r)|$ is the volume of the ball. 
The classical Dini condition, i.e., the estimate
$$
|f(x)-f(y)|\le \omega(|x-y|)
$$
with  a  continuous increasing function $\omega$ on $[0, +\infty)$ such that $\omega(0)=0$ and \eqref{din1} holds,
implies the Dini mean oscillation condition. However, the latter is weaker.
Let $f$ be a Lipschitz function outside the ball $B(0, 1/2)$, $f(0)=0$ and
$f(x)=|\ln|x||^{-\gamma}$ if $0<|x|\le 1/2$, where $0<\gamma<1$.
It is clear that the function $f$ does not satisfy the classical Dini condition. By the Poincar\'e inequality
we have $\omega(r)\le C|\ln r|^{-\gamma-1}$ with some constant $C>0$. Hence $f$ satisfies the Dini mean oscillation condition.
Obviously, a H\"older continuous function satisfies Dini's condition.

Note that any measurable function $f$ satisfying the Dini mean
oscillation condition has a version uniformly continuous on $\mathbb{R}^d$, moreover, its modulus
of continuity is estimated by $\displaystyle\int_0^{t}\frac{\omega(s)}{s}\,ds$.
We always  work with this continuous version. A survey of recent results employing the Dini mean oscillation condition is given in~\cite{Kim25}.

Let $A=(a^{ij})$ be a mapping from $\mathbb{R}^d$ to the space of symmetric matrices and let  $b=(b^i)$ be a
vector field on $\mathbb{R}^d$. We shall say that the  pair $(A, b)$ satisfies Condition~(H) if

${\rm (H_a)}$ for some $\lambda>0$
and all $x\in\mathbb{R}^d$ we have $\lambda\cdot I\le A(x)\le \lambda^{-1}\cdot I$, where $I$ is the identity matrix, and
for all $i, j\le d$ the functions $a^{ij}$ are continuous on $\mathbb{R}^d$ and satisfy the Dini mean
oscillation condition with a modulus $\omega_{a^{ij}}$.

${\rm (H_b)}$ the functions $b^i$ are Borel measurable and there exist  numbers $\beta\ge 1$, $\beta_1>0$, $\beta_2>0$, and $\beta_3>0$ such that
$$
\langle b(x), x\rangle\le \beta_1-\beta_2|x|^2, \quad |b(x)|\le \beta_3\bigl(1+|x|\bigr)^{\beta}.
$$

Note that condition ${\rm (H_a)}$ is fulfilled if for some $\lambda>0$
and all $x\in\mathbb{R}^d$ we have $\lambda\cdot I\le A(x)\le \lambda^{-1}\cdot I$ and the elements
$a^{ij}$ of the matrix $A$ have a modulus of continuity
satisfying the classical Dini condition~(\ref{din1}).

Further, condition ${\rm (H_b)}$ is fulfilled for $b(x)=-|x|^{\beta-1}x+v(x)$, where $v$ is Borel vector field and
for some numbers $C>0$, $\kappa<\beta$ we have $|v(x)|\le C\bigl(1+|x|\bigr)^{\kappa}$ for every $x\in\mathbb{R}^d$.

Set
$$
L_{A, b}u={\rm tr}(AD^2u)+\langle b, \nabla u\rangle=\sum_{i,j\le d} a^{ij}\partial_{x_i}\partial_{x_j} u+\sum_{i\le d} b^{i}\partial_{x_i} u.
$$
A function $\varrho\in L^1(\mathbb{R}^d)$ is called a probability solution to the   Kolmogorov equation
$$
\partial_{x_i}\partial_{x_j}\bigl(a^{ij}\varrho\bigr)-\partial_{x_i}\bigl(b^i\varrho\bigr)=0,
$$
which is written as
$$
L_{A, b}^{*}\varrho=0,
$$
 if
$$
\varrho\ge 0, \quad \int_{\mathbb{R}^d} \varrho(x)\,dx=1
$$
and for all $\varphi\in C_0^{\infty}(\mathbb{R}^d)$
we have
$$
\int_{\mathbb{R}^d} \varrho(x)L_{A, b}\varphi(x)\,dx=0.
$$
Note that a general Kolmogorov equation is an equation with respect to measures,
but in our situation of a non-singular diffusion matrix all solutions in the class of measures
have densities, so we consider this equation with respect to densities (in the general terminology
the measure $\varrho\, dx$ is a solution).

If $(A, b)$ satisfies Condition (H), then by virtue of \cite[Corollary 4.8]{BRSH23}
there exists a unique probability solution $\varrho$ to the equation $L_{A, b}^{*}\varrho=0$.
According to \cite[Corollary 2.3.4]{BKRS}
for each $k\ge 1$ there exists a  number $C_k$, depending only on $d$, $k$ and the quantities from Condition {\rm (H)} 
such that
$$
\int_{\mathbb{R}^d} |x|^k\varrho(x)\,dx\le C_k.
$$

\section{Main results}

We now present our main results and start with the Poisson equation. We need the following auxiliary assertion.

\begin{lemma}\label{lem1}
Suppose that a continuous function $f$ on $\mathbb{R}^d$ satisfies the Dini mean oscillation condition
with a modulus $\omega$. Let $g\in C_0^{\infty}(\mathbb{R}^d)$, $g\ge 0$, $\|g\|_{L^1(\mathbb{R}^d)}=1$.
For every $\varepsilon>0$ denote by $g_{\varepsilon}$ the function $\varepsilon^{-d}g(x/\varepsilon)$
and denote by $f*g_{\varepsilon}$ the convolution of  $f$ and $g_{\varepsilon}$.
Then the functions $f*g_{\varepsilon}$  converge uniformly to  $f$ as $\varepsilon\to 0$
and  satisfy the Dini mean oscillation condition
with the modulus~$\omega$.
\end{lemma}
\begin{proof}
Note that for every ball $B(x, r)\subset\mathbb{R}^d$ and every $y\in\mathbb{R}^d$ we have
$$
f*g_{\varepsilon}(y)-\frac{1}{|B(x, r)|}\int_{B(x, r)}f*g_{\varepsilon}(z)\,dz=
\int_{\mathbb{R}^d}\bigl(f(y-z)-f_B(x-z, r)\bigr)g_{\varepsilon}(z)\,dz.
$$
Then
\begin{multline*}
\frac{1}{|B(x, r)|}\int_{B(x, r)}\biggl|f*g_{\varepsilon}(y)-\frac{1}{|B(x, r)|}\int_{B(x, r)}f*g_{\varepsilon}(z)\,dz\biggr|\,dy
\\
\le\int_{\mathbb{R}^d}\biggl(\frac{1}{|B(x, r)|}\int_{B(x, r)}\bigl|f(y-z)-f_B(x-z, r)\bigr|\,dy\biggr) g_{\varepsilon}(z)\,dz,
\end{multline*}
where the right-hand side is estimated by $\omega(r)$.

The functions $f*g_{\varepsilon}$ converge uniformly to $f$ because
$f$ is a continuous function and these functions have a common modulus of continuity.
\end{proof}

The following theorem is our first main result.
For a domain $\Omega$ in $\mathbb{R}^d$ we denote by $W^{p,2}(\Omega)$ 
the classical Sobolev space of functions $f$ from $L^p(\Omega)$ 
having first and second Sobolev derivatives $\partial_{x_i}f$ and $\partial_{x_i}\partial_{x_j}f$ in $L^p(\Omega)$. 
It is equipped with the  Sobolev norm 
$$
\|f\|_{W^{p,2}(\Omega)}=\|f\|_{L^{p}(\Omega)}+\sum_{i\le d}+\|\partial_{x_i}f\|_{L^{p}(\Omega)}+
\sum_{i,j\le d} \|\partial_{x_i}\partial_{x_j}\|_{L^{p}(\Omega)}.
$$
The class $W^{p,2}_{loc}(\mathbb{R}^d)$ of locally Sobolev functions consists of all functions $f$ such 
that $\zeta f\in W^{p,2}(\mathbb{R}^d)$ for all $\in C_0^\infty(\mathbb{R}^d)$. 

\begin{theorem}\label{th1}
Suppose that $(A, b)$ satisfies  Condition {\rm (H)} and $\varrho$ is a probability solution to the equation $L_{A, b}^{*}\varrho=0$.
Let $p>d$, $k\ge 1$ and let  $\psi\in C_b(\mathbb{R}^d)$ be such that
 $$
 \int_{\mathbb{R}^d} \psi(x)\varrho(x)\,dx=0.
 $$
Then there exists a function $u\in W^{p,2}_{loc}(\mathbb{R}^d)\bigcap C^1(\mathbb{R}^d)$ such that
$$L_{A, b}u=\psi$$
 and
$$
|u(x)|\le C(1+|x|^k)\sup_{y\in\mathbb{R}^d}\frac{|\psi(y)|}{1+|y|^k}, \quad
|\nabla u(x)|\le C(1+|x|^{k+\beta})\sup_{y\in\mathbb{R}^d}\frac{|\psi(y)|}{1+|y|^k},
$$
$$
\biggl(\int_{\mathbb{R}^d} \frac{|D^2u(x)|^p}{1+|x|^{s}}\,dx\biggr)^{1/p}\le C\sup_{y\in\mathbb{R}^d}\frac{|\psi(y)|}{1+|y|^k}.
$$
where $C$ and $s$ depend only on $p$, $d$, and the quantities from Condition {\rm (H)}.
\end{theorem}
\begin{proof}
Suppose first that the functions $a^{ij}$
are continuously differentiable.
According to \cite[Theorem 4.1 and Example 4.2 (i)]{BRSH18}, there is a solution
$u\in W^{p, 2}_{loc}(\mathbb{R}^d)\bigcap C^1(\mathbb{R}^d)$ to the equation $L_{A, b}u=\psi$ such that
this solution is growing at infinity not faster than $\ln(1+|x|)$. Moreover, by \cite[Proposition~4.1]{BRSH18}
the function $|\nabla u|$ has an at most polynomial growth at infinity.
Note that the  constants in the corresponding estimates  obtained in \cite{BRSH18}
depend not only on the quantities from Condition~(H), but  also on the Sobolev norms of $a^{ij}$.
Our further reasoning aims at obtaining estimates not involving Sobolev norms of~$a^{ij}$.

Let $\varphi\in C_0^{\infty}(\mathbb{R}^d)$, $\varphi\ge 0$ and $\varphi(x)=1$ if $|x|<1$.
Set $\varphi_N(x)=\varphi(x/N)$. The following equality holds:
$$
L_{A, b}(u^2\varphi_N)-u^2L_{A, b}\varphi_N-4u\langle A\nabla u, \nabla\varphi_N\rangle=2|\sqrt{A}\nabla u|^2\varphi_N+2u\psi\varphi_N.
$$
Multiplying this equality by $\varrho$ and integrating, we obtain
$$
\int_{\mathbb{R}^d} \bigl(-u^2L_{A, b}\varphi_N-4u\langle A\nabla u, \nabla\varphi_N\rangle\bigr)\varrho\,dx=
\int_{\mathbb{R}^d} \bigl(2|\sqrt{A}\nabla u|^2\varphi_N+2u\psi\varphi_N\bigr)\varrho\,dx.
$$
Note that throughout we omit indication of variables of functions under the integral sign in case of long expressions. 
Observe that the functions
$$
|a^{ij}|u^2, \quad u^2|b^i|, \quad |A\nabla u|, \quad |\sqrt{A}\nabla u|^2, \quad u\psi
$$
are integrable with respect to the  measure $\varrho\,dx$ on $\mathbb{R}^d$.
Therefore, letting $N\to\infty$, we arrive at  the equality
$$
\int_{\mathbb{R}^d} 2|\sqrt{A}\nabla u|^2\varrho\,dx=-2\int_{\mathbb{R}^d} u\psi\varrho\,dx.
$$
According to \cite[Theorem 3.5]{BRSH23}, for every ball $B(0, R)$ there exists
a positive number $Q_R$, depending only on the quantities from Condition~(H), such that the following Harnack  inequality
is fulfilled:
$$
\sup_{B(0, R)}\varrho(x)\le Q_R\inf_{B(0, R)}\varrho(x).
$$
Since
$$
L_{A, b}(|x|^2)\le 2d\lambda+2\beta_1-2\beta_2|x|^2,
$$
there  exist  numbers $R_0>0$ and $M_0>0$ such that
$$
L_{A, b}(|x|^{2k})\le -M_0(1+|x|^{2k})\le -M_0\quad \hbox{if \ } |x|>R_0.
$$
The numbers $M_0$ and $R_0$ depend only on the numbers $d$, $k$, $\lambda$, $\beta_1$, and $\beta_2$.
Then by \cite[Corollary~2.3.3]{BKRS} we have
$$
\int_{|x|>2R_0}\varrho(x)\,dx\le M_0^{-1}\int_{B(0, 2R_0)}\varrho(x)\bigl|L_{A, b}(|x|^{2k})\bigr|\,dx,
$$
 which yields the estimate
$$
1\le |B(0, 2R_0)|\Bigl(1+M_0^{-1}\sup_{B(0, 2R_0)}\bigl|L_{A, b}(|x|^{2k})\bigr|\Bigr)\sup_{B(0, 2R_0)}\varrho(x).
$$
Combining this estimate and Harnack's inequality we obtain
$$
\inf_{B(0, 2R_0)}\varrho(x)\ge C_1>0,
$$
where the number $C_1$ depends only on the quantities from Condition (H).
Therefore, one has the estimate
$$
\int_{B(0, 2R_0)}|\nabla u|^2\,dx\le C_1^{-1}\lambda^{-1}\int_{\mathbb{R}^d} |u||\psi|\varrho\,dx.
$$
Subtracting from the function $u$ its average over the ball $B(0, 2R_0)$, we assume further  that
$$
\int_{B(0, 2R_0)}u(x)\,dx=0.
$$
By the Poincar\'e  inequality
$$
\int_{B(0, 2R_0)}|u|^2\,dx\le C_2\int_{\mathbb{R}^d} |u||\psi|\varrho\,dx,
$$
where $C_2$ depends only on $R_0$, $d$ and $C_1$.
Estimating $\sup_{B(0, R_0)}|u|$ by the $L^2$-norm of $u$ on the ball $B(0, 2R_0)$ according to
\cite[Corollary 9.21]{GTrud}, we obtain
$$
\sup_{B(0, R_0)}|u(x)|\le C_3\biggl(\int_{\mathbb{R}^d}|\psi||u|\varrho\,dx\biggr)^{1/2}+C_3\sup_{B(0, R_0)}|\psi(x)|,
$$
where $C_3$ depends only on $d$, $k$ and the quantities from Condition (H).

Set $u^2=wV$, where $V(x)=1+|x|^{2k}$. Then
$$
L_{A, b}w+2V^{-1}\langle A\nabla V, \nabla w\rangle+w\frac{L_{A, b}V}{V}\ge 2u\psi V^{-1}.
$$
We recall that $L_{A, b}V(x)\le -M_0V(x)$ whenever $|x|>R_0$. Therefore,
$$
L_{A, b}w+2V^{-1}\langle A\nabla V, \nabla w\rangle-M_0 w\ge 2u\psi V^{-1}-C_4\sup_{B(0, R_0)}|u(x)|^2,
$$
where
$$
C_4=M_0+\sup_{B(0, R_0)}\frac{|L_{A, b}V(x)|}{V(x)}.
$$
Since the function $w$ tends to zero at infinity, it has a  point $y$ of global maximum.
Applying \cite[Theorem 9.6]{GTrud}, we obtain
$$
-M_0w(y)\ge -2|u(y)||\psi(y)|V(y)^{-1}-C_4\sup_{B(0, R_0)}|u(x)|^2.
$$
Therefore, we have
\begin{multline*}
\sup_{\mathbb{R}^d}\frac{u^2(x)}{V(x)}\le 2M_0^{-1}
\sup_{\mathbb{R}^d}\frac{|u(x)|}{\sqrt{V(x)}}\cdot\sup_{\mathbb{R}^d} \frac{|\psi(x)|}{\sqrt{V(x)}}
\\
+2C_4C_3^2M_0^{-1}\int_{\mathbb{R}^d}|\psi(x)||u(x)|\varrho(x)\,dx+2C_4C_3^2\sup_{B(0, R_0)}|\psi(x)|^2.
\end{multline*}
Since
$$
\int_{\mathbb{R}^d}|\psi(x)||u(x)|\varrho(x)\,dx\le
\sup_{\mathbb{R}^d}\frac{|u(x)|}{\sqrt{V(x)}}\cdot
\sup_{\mathbb{R}^d}\frac{|\psi(x)|}{\sqrt{V(x)}}\int_{\mathbb{R}^d}V(x)\varrho(x)\,dx,
$$
 we arrive at the estimate
$$
\sup_{\mathbb{R}^d}\frac{|u(x)|}{\sqrt{1+|x|^{2k}}}\le C_5\sup_{\mathbb{R}^d}\frac{|\psi(x)|}{\sqrt{1+|x|^{2k}}}.
$$
Note that $C_5$ depends only on $d$, $k$ and the quantities from our conditions on the coefficients.

Let us estimate the first derivatives of the solution $u$.
Set
$$
q=1+\sup_{B(z, 4)}|b(x)|.
$$
Let $u(x)=v(q(x-t))$ and $v(y)=u(q^{-1}y+t)$, where $t\in B(z, 1)$.
Then $v$ solves the equation
$$
\widetilde{a}^{ij}(y)=a^{ij}(q^{-1}y+t),
\quad \widetilde{b}^i(y)=q^{-1}b^i(q^{-1}y+t), \quad \widetilde{\psi}(y)=q^{-2}\psi(q^{-1}y+t).
$$
Observe that $|\widetilde{b}(y)|\le 1$ if $y\in B(0, 2)$. By \cite[Theorem 9.11]{GTrud} we have
$$
\|v\|_{W^{p, 2}(B(0, 1))}\le C(d, p)\Bigl(\sup_{B(0, 2)}|v(y)|+\sup_{B(0, 2)}|\widetilde{\psi}(y)|\Bigr),
$$
where $C(d, p)$ depends only on $d$ and $p$. Applying the Sobolev  embedding theorem, we arrive at  the estimate
$$
\sup_{B(0, 2)}|\nabla v(y)|\le \widetilde{C}(d, p)\Bigl(\sup_{B(0, 2)}|v(y)|+\sup_{B(0, 2)}|\widetilde{\psi}(y)|\Bigr).
$$
Returning to the original coordinates, we obtain
$$
|\nabla u(t)|\le \widetilde{C}(d, p)\Bigl(q\sup_{B(t, 2)}|u(x)|+q^{-1}\sup_{B(t, 2)}|\psi(x)|\Bigr).
$$
In particular, this estimate is fulfilled if $t=z$.
Taking into account the inequalities obtained above and the condition
$$|b(x)|\le \beta_3(1+|x|)^{\beta},$$
we arrive at the estimate
$$
|\nabla u(z)|\le C_6(1+|z|)^{\beta+k}\sup_{\mathbb{R}^d}\frac{|\psi(x)|}{1+|x|^k}.
$$
We now estimate the second derivatives of the solution. From the aforementioned estimate 
of the Sobolev norm we have 
$$
\int_{B(0, 1)}|D^2v(y)|^p\,dy\le C(d, p)^p\Bigl(\sup_{B(0, 2)}|v(y)|+\sup_{B(0, 2)}|\widetilde{\psi}(y)|\Bigr)^p.
$$
Returning to the original coordinates, this gives for every $t\in B(z, 1)$
the inequality
$$
q^d\int_{B(t, q^{-1})}q^{-2p}|D^2u(x)|^p\,dx\le
C(d, p)^p\Bigl(\sup_{B(t, 2)}|u(x)|+q^{-2}\sup_{B(t, 2)}|\psi(x)|\Bigr)^p.
$$
Therefore, we have the estimate
$$
\int_{B(z, 1)}|D^2u(x)|^p\,dx\le
C_7(1+|z|)^{(2\beta+k)p}\Biggl(\sup_{\mathbb{R}^d}\frac{|\psi(x)|}{1+|x|^k}\Biggr)^p,
$$
which yields the bound
$$
\int_{B(z, 1)}\frac{|D^2u(x)|^p}{1+|x|^s}\,dx\le
\frac{C_8}{(1+|z|)^{d+1}}\Biggl(\sup_{\mathbb{R}^d}\frac{|\psi(x)|}{1+|x|^k}\Biggr)^p,
\quad s=(2\beta+k)p+d+1.
$$
Integrating in $z$ and applying Fubini's theorem, we obtain
$$
\int_{\mathbb{R}^d}\frac{|D^2u(x)|^p}{1+|x|^s}\,dx\le
C_9\Biggl(\sup_{\mathbb{R}^d}\frac{|\psi(x)|}{1+|x|^k}\Biggr)^p.
$$
Thus, in the case of continuously differentiable functions $a^{ij}$ we have constructed a solution
for which the desired estimates are fulfilled with constants depending only on $d$, $k$, $p$ and the quantities from Condition~(H).

We now  consider the case where $(A, b)$ satisfies  Condition (H) and there are no assumptions about additional smoothness of the
coefficients. By Lemma \ref{lem1} there exists a sequence of matrices $\{A_n\}$ such that $\{A_n\}$
converges uniformly to $A$ and for every $n$ the matrix
$A_n$ satisfies Condition ${\rm (H_a)}$ and its elements $a^{ij}_n$ are continuously differentiable functions.
As shown above, the equation $L_{A_n, b}u=\psi$ has a solution $u_n$ satisfying the estimate from the  assertion of the theorem.
Then the sequence $\{u_n\}$ contains a  subsequence $\{u_{n_j}\}$ such that on every ball $B(0, m)$
the sequence $\{u_{n_j}\}$ converges to $C^1(B(0, m))$ and weakly in $W^{p, 2}(B(0, m))$
to some function $u\in W^{p,2}_{loc}(\mathbb{R}^d)\bigcap C^1(\mathbb{R}^d)$. 
It is clear that $u$ is a solution to the equation $L_{A, b}u=\psi$
and the   estimates indicated above are valid for~$u$.
\end{proof}

Let us apply the previous theorem  for estimating distances between solutions to stationary  Kolmogorov equations.

We need the following auxiliary assertion. Recall that
for each $k\ge 1$ there exists a  number $C_k$ such that
$$
\int_{\mathbb{R}^d} |x|^k\varrho(x)\,dx\le C_k.
$$

\begin{lemma}\label{lem2}
Suppose that $\varrho$ is a probability solution to the equation $L_{A, b}^{*}\varrho=0$ and Condition~{\rm (H)}
is fulfilled. Then, for every $k\ge 1$, there exists a number $M_k$, depending only on $k$, $\lambda$, $d$,
$\beta$, $\beta_1$, $\beta_2$, and $\beta_3$ such that
$$
\|(1+|x|)^k\varrho\|_{L^{d/(d-1)}(\mathbb{R}^d)}\le M_k.
$$
\end{lemma}
\begin{proof}
By \cite[Corollary 1.5.3 and Remark 1.5.5 (ii)]{BKRS} for every $z\in\mathbb{R}^d$
the inequality
$$
\int_{B(z, 1)}\varrho(x)^{d/(d-1)}\,dx\le C_1\Bigl(\int_{B(z, 2)}\bigl(1+|b(x)|\bigr)\varrho(x)\,dx\Bigr)^{d/(d-1)}
$$
holds with a number $C_1>0$ depending only on $\lambda$ and $d$. Applying H\"older's inequality to the right-hand side we
obtain
$$
\int_{B(z, 1)}\varrho(x)^{d/(d-1)}\,dx\le C_1\int_{B(z, 2)}\bigl(1+|b(x)|\bigr)^{d/(d-1)}\varrho(x)\,dx.
$$
Multiplying this inequality by $(1+|z|)^{dk/(d-1)}$ and using the estimates
$$
\frac{1}{3}\le \inf_{x\in B(z, 2)}\frac{1+|z|}{1+|x|}\le \sup_{x\in B(z, 2)}\frac{1+|z|}{1+|x|}\le 3,
$$
we arrive at the inequality
\begin{multline*}
\int_{B(z, 1)}(1+|x|)^{dk/(d-1)}\varrho(x)^{d/(d-1)}\,dx
\\
\le 9^{dk/(d-1)}C_1\int_{B(z, 2)}\bigl(1+|x|\bigr)^{dk/(d-1)}\bigl(1+|b(x)|\bigr)^{d/(d-1)}\varrho(x)\,dx.
\end{multline*}
According to  condition $\rm(H_b)$ we have
$$
\bigl(1+|x|\bigr)^{dk/(d-1)}\bigl(1+|b(x)|\bigr)^{d/(d-1)}\le C_2(1+|x|)^{(dk+\beta k)/(d-1)},
$$
where $C_2$ depends only on $k$, $d$, $\beta$, $\beta_3$.
Hence
$$
\int_{B(z, 1)}(1+|x|)^{dk/(d-1)}\varrho(x)^{d/(d-1)}\,dx\le
9^{dk/(d-1)}C_1C_2\int_{B(z, 2)}(1+|x|)^{(dk+\beta k)/(d-1)}\varrho(x)\,dx.
$$
Integrating with respect to $z$ and applying Fubini's theorem, we obtain
\begin{multline*}
|B(0, 1)|\int_{\mathbb{R}^d}(1+|x|)^{dk/(d-1)}\varrho(x)^{d/(d-1)}\,dx
\\
\le 9^{dk/(d-1)}C_1C_2|B(0, 2)|\int_{\mathbb{R}^d}(1+|x|)^{(dk+\beta k)/(d-1)}\varrho(x)\,dx,
\end{multline*}
where the right-hand side is estimated by a number depending  only on $k$, $\lambda$, $d$ and
the quantities from $\rm (H_b)$.
\end{proof}

The next theorem is the second main result of the paper.

\begin{theorem}\label{th2}
Suppose that  $(A_{\mu}, b_{\mu})$ and $(A_{\sigma}, b_{\sigma})$ satisfy  Condition~{\rm (H)}.
Let $\varrho_{\mu}$ and $\varrho_{\sigma}$ be probability solutions to the equations
$L_{A_{\mu}, b_{\mu}}^{*}\varrho_{\mu}=0$ and $L_{A_{\sigma}, b_{\sigma}}^{*}\varrho_{\sigma}=0$.
Then, for every $r>1$ and $k\ge 1$, we have the estimate
\begin{multline*}
\|(1+|x|^k)(\varrho_{\mu}-\varrho_{\sigma})\|_{L^1(\mathbb{R}^d)}\le
C\|A_{\mu}-A_{\sigma}\|_{L^r(\mathbb{R}^d, \varrho_{\sigma}\,dx)}
\\
+
C\int_{\mathbb{R}^d}|b_{\mu}(x)-b_{\sigma}(x)|(1+|x|^{\beta+k})\varrho_{\sigma}(x)\,dx,
\end{multline*}
where the number $C$ depends only on $d$, $r$, $k$, and the quantities from Condition {\rm (H)}.
\end{theorem}
\begin{proof}
Since $\varrho_{\sigma}$ is a probability density, it suffices to obtain this estimate  for all numbers $r$ in
$(1,d/(d-1))$.

Observe that for every function $v\in C_0^{\infty}(\mathbb{R}^d)$ we have the equality
$$
\int_{\mathbb{R}^d} \bigl(\varrho_{\mu}-\varrho_{\sigma}\bigr)L_{A_{\mu}, b_{\mu}}v\,dx=
-\int_{\mathbb{R}^d} \bigl[{\rm tr}\bigl((A_{\mu}-A_{\sigma})D^2v\bigr)+
\langle b_{\mu}-b_{\sigma}, \nabla v\rangle\bigr]\varrho_{\sigma}\,dx.
$$
This equality remains true for all functions $v\in W^{p,2}_{loc}(\mathbb{R}^d)\bigcap C^1(\mathbb{R}^d)$, $p>d$,
with compact support.

Let $\varphi\in C_0^{\infty}(\mathbb{R}^d)$ and $|\varphi(x)|(1+|x|)^{-k}\le 1$.
We denote by $u\in W^{p,2}_{loc}(\mathbb{R}^d)\bigcap C^1(\mathbb{R}^d)$ the solution to the  Poisson equation
$$L_{A_{\mu}, b_{\mu}}u=\psi,$$
constructed in Theorem \ref{th1}, where
$$
\psi(x)=\varphi(x)-\int_{\mathbb{R}^d} \varphi(x)\varrho_{\mu}(x)\,dx.
$$
The number $p>d$ will be chosen below.
Note that
$$
\frac{|\psi(x)|}{1+|x|^k}\le 1+\int_{\mathbb{R}^d} (1+|x|^k)\varrho_{\mu}(x)\,dx\le C_1,
$$
where $C_1$ depends only on the quantities from Condition (H).

Let $\zeta\in C_0^{\infty}(\mathbb{R}^d)$, $\zeta\ge 0$ and $\zeta(x)=1$ if $|x|<1$.
Set $\zeta_N(x)=\zeta(x/N)$. Then
$$
\int_{\mathbb{R}^d} \bigl(\varrho_{\mu}-\varrho_{\sigma}\bigr)L_{A_{\mu}, b_{\mu}}(\zeta_Nu)\,dx=
-\int_{\mathbb{R}^d} \bigl[{\rm tr}\bigl((A_{\mu}-A_{\sigma})D^2(\zeta_Nu)\bigr)+
\langle b_{\mu}-b_{\sigma}, \nabla(\zeta_Nu)\rangle\bigr]\varrho_{\sigma}\,dx.
$$
Observe that the functions
$$
{\rm tr}(A_{\mu}D^2u), \quad {\rm tr}(A_{\sigma}D^2u), \quad |A_{\mu}\nabla u|, \quad |A_{\sigma}\nabla u|
\quad \langle b_{\mu}, \nabla u\rangle, \quad \langle b_{\sigma}, \nabla u\rangle
$$
are integrable with respect to $\varrho_{\mu}\,dx$ and $\varrho_{\sigma}\,dx$ on $\mathbb{R}^d$.
Letting $N\to\infty$, we arrive at the equality
$$
\int_{\mathbb{R}^d} \psi(\varrho_{\mu}-\varrho_{\sigma})\,dx=
-\int_{\mathbb{R}^d} \bigl[{\rm tr}\bigl((A_{\mu}-A_{\sigma})D^2u\bigr)+
\langle b_{\mu}-b_{\sigma}, \nabla u\rangle\bigr]\varrho_{\sigma}\,dx.
$$
Since $\varrho_{\mu}$ and $\varrho_{\sigma}$ are probability densities, we have
$$
\int_{\mathbb{R}^d} \psi(\varrho_{\mu}-\varrho_{\sigma})\,dx=\int_{\mathbb{R}^d} \varphi(\varrho_{\mu}-\varrho_{\sigma})\,dx.
$$
By H\"older's  inequality
$$
-\int_{\mathbb{R}^d} {\rm tr}\bigl((A_{\mu}-A_{\sigma})D^2u\bigr)\varrho_{\sigma}\,dx\le \|(A_{\mu}-A_{\sigma})\|_{L^r(\mathbb{R}^d, \varrho_{\sigma}\,dx)}
\|D^2u\|_{L^{r'}(\mathbb{R}^d, \varrho_{\sigma}\,dx)},
$$
where
$$
 r'=\frac{r}{r-1}>d.
$$
Let $p=r'd$ and let $s$ be the number from the Theorem~\ref{th1}. We have
$$
\int_{\mathbb{R}^d}|D^2u(x)|^{r'}\varrho_{\sigma}(x)\,dx=
\int_{\mathbb{R}^d}\frac{|D^2u(x)|^{r'}}{(1+|x|)^{s/d}}(1+|x|)^{s/d}\varrho_{\sigma}(x)\,dx.
$$
Applying again H\"older's inequality we obtain
\begin{multline*}
\int_{\mathbb{R}^d}|D^2u(x)|^{r'}\varrho_{\sigma}(x)\,dx
\\
\le\biggl(\int_{\mathbb{R}^d}\frac{|D^2u(x)|^{p}}{(1+|x|)^{s}}\,dx\biggr)^{1/d}
\biggl(\int_{\mathbb{R}^d}(1+|x|)^{s/(d-1)}\varrho_{\sigma}(x)^{d/(d-1)}\,dx\biggr)^{(d-1)/d}.
\end{multline*}
According to Lemma~\ref{lem2},
we have
$$
\|D^2u\|_{L^{r'}(\mathbb{R}^d, \varrho_{\sigma}\,dx)}^{r'}\le C_2,
$$
where $C_2$ depends only on $r$, $d$ and the quantities from Condition (H).

Since  $|\nabla u|$ satisfies the inequality
$$
|\nabla u(x)|\le C_3(1+|x|^{\beta+k}),
$$
where $C_3$ depends only on $r$, $d$ and the quantities from Condition (H), we have
$$
-\int_{\mathbb{R}^d} \langle b_{\mu}-b_{\sigma}, \nabla u\rangle\varrho_{\sigma}\,dx\le
C_3\int_{\mathbb{R}^d} |b_{\mu}-b_{\sigma}|(1+|x|^{\beta+k})\varrho_{\sigma}\,dx.
$$
Combining the obtained estimates, we arrive at the inequality
$$
\int_{\mathbb{R}^d} \varphi(\varrho_{\mu}-\varrho_{\sigma})\,dx\le
C_2\|(A_{\mu}-A_{\sigma})\|_{L^r(\mathbb{R}^d, \varrho_{\sigma}\,dx)}
+C_3\int_{\mathbb{R}^d} |b_{\mu}-b_{\sigma}|(1+|x|^{\beta+k})\varrho_{\sigma}\,dx,
$$
from which  the assertion of theorem follows, because $\varphi$ was arbitrary.
\end{proof}

\section{An application to nonlinear equations}

Here we consider an application of the last  theorem to nonlinear Kolmogorov equations (see \cite{BS24} for a recent survey).

Let $k\ge 1$. We denote by $\mathcal{P}_k(\mathbb{R}^d)$ the space of
probability densities $\varrho$ on $\mathbb{R}^d$ such that $|x|^k\varrho\in L^1(\mathbb{R}^d)$.
This is a complete metric space with respect to the metric
$$\|(1+|x|^k)(\varrho_1-\varrho_2)\|_{L^1(\mathbb{R}^d)}.$$

Suppose that for every $\varepsilon\in[0, 1]$ we are given functions
$$
(x, \varrho)\mapsto a^{ij}_{\varepsilon}(x, \varrho), \quad
(x, \varrho)\mapsto b^i_{\varepsilon}
$$
on the space $\mathbb{R}^d\times\mathcal{P}_k(\mathbb{R}^d)$. Assume that
these functions are Borel measurable with respect to the variable~$x$ and the matrices
$A_{\varepsilon}(x, \varrho)=\bigl(a^{ij}_{\varepsilon}(x, \varrho)\bigr)$ are symmetric and nonnegative definite.

Let us consider the nonlinear Kolmogorov equation
$$
\partial_{x_i}\partial_{x_j}\bigl(a^{ij}_{\varepsilon}(x, \varrho)\varrho\bigr)
-\partial_{x_i}\bigl(b^i_{\varepsilon}(x, \varrho)\varrho\bigr)=0.
$$
Using the operator
$$
L_{\varepsilon, \varrho}u(x)={\rm tr}(A_{\varepsilon}(x, \varrho)D^2u(x))+\langle b_{\varepsilon}(x, \varrho), \nabla u(x)\rangle
$$
we write the equation in the short form
\begin{equation}\label{nonl}
L_{\varepsilon, \varrho}^{*}\varrho=0.
\end{equation}
We say that $\varrho\in\mathcal{P}_k(\mathbb{R}^d)$ is a solution to equation~(\ref{nonl})
if for all $\varphi\in C_0^{\infty}(\mathbb{R}^d)$
we have
$$
\int_{\mathbb{R}^d} \varrho(x)L_{\varepsilon, \varrho}\varphi(x)\,dx=0,
$$
which means that $\varrho$ is a solution to the linear Kolmogorov equation with the operator $L_{\varepsilon, \varrho}$
depending on $\varrho$.

\begin{theorem}\label{th3}
Suppose that for all $\varepsilon\in[0, 1]$ and $\varrho\in\mathcal{P}_k(\mathbb{R}^d)$
the functions
$$
x\mapsto a^{ij}_{\varepsilon}(x, \varrho)\quad {\rm and} \quad x\mapsto b^i(x, \varrho)
$$
satisfy Condition~{\rm(H)} and the quantities $\lambda$, $w(r)$, $\beta$, $\beta_1$, $\beta_2$, $\beta_3$
do not depend on $\varepsilon$ and $\varrho$. Assume also that there exist $N>0$ and $m\ge 0$ such that
the estimate
$$
|A(x, \varrho_1)-A(x, \varrho_2)|+|b(x, \varrho_1)-b(x, \varrho_2)|\le \varepsilon N(1+|x|^m)\|(1+|x|^k)(\varrho_1-\varrho_2)\|_{L^1(\mathbb{R}^d)}
$$
is fulfilled for all $\varepsilon\in[0, 1]$, $x\in\mathbb{R}^d$ and $\varrho_1, \varrho_2\in\mathcal{P}_k(\mathbb{R}^d)$.
Then there exists a positive number $\varepsilon_0$
such that for every $\varepsilon\in[0, \varepsilon_0]$ there exists a unique solution $\varrho\in \mathcal{P}_k(\mathbb{R}^d)$
to  equation~{\rm(\ref{nonl})}.
\end{theorem}
\begin{proof}
As noted in Section~2,  for all $\varepsilon\in[0, 1]$ and $\varrho\in\mathcal{P}_k(\mathbb{R}^d)$ the linear
Kolmogorov equation $L_{\varepsilon, \varrho}^{*}\sigma=0$ has a unique solution $\sigma\in\mathcal{P}_k(\mathbb{R}^d)$.
Let us consider the mapping
$$\Phi_{\varepsilon}\colon \mathcal{P}_k(\mathbb{R}^d)\to \mathcal{P}_k(\mathbb{R}^d)$$
such that
$$
\sigma=\Phi_{\varepsilon}(\varrho) \quad \Longleftrightarrow \quad L_{\varepsilon, \varrho}^{*}\sigma=0.
$$
There exists a number $M>0$  depending only on $k$, $m$, $d$ and the quantities from Condition~(H)
such that for all $\varepsilon\in[0, 1]$ and $\varrho\in\mathcal{P}_k(\mathbb{R}^d)$ we have
$$
\int_{\mathbb{R}^d}(1+|x|)^{2m+\beta+k}\Phi_{\varepsilon}(\varrho)(x)\,dx\le M.
$$
By Theorem~\ref{th2} with $r=2$ we have the estimate
\begin{multline*}
\|(1+|x|^k)(\Phi_{\varepsilon}(\varrho_1)-\Phi_{\varepsilon}(\varrho_2))\|_{L^1(\mathbb{R}^d)}\le
C\|A(\,\cdot\,, \varrho_1)-A(\,\cdot\,, \varrho_2)\|_{L^2(\mathbb{R}^d, \Phi_{\varepsilon}(\varrho_2)\,dx)}
\\
+
C\int_{\mathbb{R}^d}|b(x, \varrho_1)-b(x, \varrho_2)|(1+|x|^{\beta+k})\Phi_{\varepsilon}(\varrho_2)(x)\,dx.
\end{multline*}
According to our assumptions we obtain
$$
\|(1+|x|^k)(\Phi_{\varepsilon}(\varrho_1)-\Phi_{\varepsilon}(\varrho_2))\|_{L^1(\mathbb{R}^d)}\le
\varepsilon NC(\sqrt{M}+M)\|(1+|x|^k)(\varrho_1-\varrho_2)\|_{L^1(\mathbb{R}^d)},
$$
Let $\varepsilon_0>0$ be such that $\varepsilon_0NC(\sqrt{M}+M)<1$. Then for every $\varepsilon\in[0, \varepsilon_0]$
the mapping $\Phi_{\varepsilon}$ is a contraction.
According to the Banach theorem there exists a unique element $\varrho\in\mathcal{P}_k(\mathbb{R}^d)$
such that $\varrho=\Phi_{\varepsilon}(\varrho)$.
\end{proof}

Let us consider an example where all hypotheses of Theorem~\ref{th3} are fulfilled.

\begin{example}
\rm
Let $A_0=(a^{ij}_0)$ be a mapping from $\mathbb{R}^d$ to the space of symmetric matrices and let  $b_0=(b_0^i)$ be a
vector field on $\mathbb{R}^d$. Suppose that the  pair $(A_0, b_0)$ satisfies Condition~(H). Set
$$
A_{\varepsilon}(x, \varrho)=A_0(x)+\varepsilon\int_{\mathbb{R}^d}q(x, y)\varrho(y)\,dy, \quad
b_{\varepsilon}(x, \varrho)=b_0(x)+\varepsilon\int_{\mathbb{R}^d}h(x, y)\varrho(y)\,dy,
$$
where $q=(q^{ij})$, the  functions $q^{ij}$ are Borel measurable, bounded and H\"older continuous~in~$x$ uniformly with respect to~$y$,
$h=(h^i)$, and $h^i$  are bounded Borel measurable functions. Then the pair $(A_{\varepsilon}, b_{\varepsilon})$
satisfies all assumptions of Theorem~\ref{th3}.
\end{example}

\end{document}